\documentclass[12pt,a4paper]{article}
\usepackage{amsthm}
\usepackage{amssymb}
\usepackage[dvips]{graphicx}

\renewcommand{\phi}{\varphi}
\renewcommand{\epsilon}{\varepsilon}

\newcommand{\grad}{\mathop{\mathrm{grad}}}
\newcommand{\tr}{\mathop{\mathrm{tr}}}

\newcommand{\R}{\mathbb{R}}

\newcommand{\pa}[1]{{\frac\partial{\partial#1}}}

\newtheorem{theorem}{Theorem}
\newtheorem{lemma}{Lemma}[section]
\newtheorem{proposition}{Proposition}[section]
\newtheorem{definition}{Definition}
\newtheorem{corollary}{Corollary}
\newtheorem{remark}{Remark}[section]

\author{Roman M. Fedorov\thanks{Partially supported by the RFBR grant 02--01--00482}}

\title{
Upper bounds\\
for the number of orbital topological types\\
of planar polynomial vector fields\\
``modulo limit cycles''}

\begin{document}
\thispagestyle{empty}

\maketitle
\begin{abstract}
The paper deals with planar polynomial vector fields. We aim to estimate the
number of orbital topological equivalence classes for the fields of degree
$n$\/. An evident obstacle for this is the second part of Hilbert's 16th
problem. To circumvent this obstacle we introduce the notion of {\em
equivalence modulo limit cycles\/}. This paper is the continuation of the
author's paper in [Mosc. Math. J. \textbf{1} (2001), no. 4] where the lower
bound of the form $2^{cn^2}$ has been obtained. Here we obtain the upper
bound of the same form. We also associate an equipped planar graph to every
planar polynomial vector field, this graph is a complete invariant for
orbital topological classification of such fields.
\end{abstract}

\tableofcontents{}
\vfill\eject

\section{Introduction}

Let us consider a planar polynomial vector field:
\[
v(x,y)=P(x,y)\pa x+Q(x,y)\pa y.
\]
Recall that vector fields~$v_1$ and~$v_2$ are called \emph{orbitally
topologically equivalent\/} if there is a homomorphism~$\phi$ of the phase
space such that~$\phi$ takes any trajectory of~$v_1$ to a trajectory of~$v_2$
and~$\phi$ preserves the natural orientation of these trajectories
(see~\cite{ds1} and~\cite{Arn2}).

The phase space of a vector field below is the plane, or a closed subspace of
the plane such that its boundary consists of limit cycles of the field (see
Definition~\ref{pfield}).

It follows from Bezout's theorem that a singular point set of the polynomial
vector field~$v$ is either finite or contains an algebraic curve. In the
first case this set consists of at most $(\deg v)^2$ points. In the second
case the equation of the algebraic curve is given by
\[
F(x,y)=0,
\]
where~$F$ is the greatest common divisor of~$P$ and~$Q$. We are mainly
interested in the first case. Notice that we can, in some sense, reduce the
second case to the first one, dividing~$P$ and~$Q$ by~$F$.

\begin{definition}\label{asterisk}\rm
A planar polynomial vector field is called \emph{$*$-field\/} if it has only
a finite number of singular points.
\end{definition}

\subsection{Equivalence modulo limit cycles}

One of the questions in the second part of the 16th Hilbert Problem is the
following: Is it true that for any~$n$ there is a constant $H(n)$ such that
every polynomial vector field of degree at most $n$ has at most $H(n)$ limit
cycles? The answer is not known even for $n=2$.

Our goal is to find some upper and lower bounds for the number of orbital
topological equivalence classes of degree $n$ planar polynomial vector
fields. The Hilbert Problem is an obstacle to obtaining such an estimate
because orbitally topologically equivalent vector fields have the same number
of limit cycles (see Corollary~\ref{hilbert}). To avoid this difficulty we
shall introduce the concept of \emph{equivalence modulo limit cycles\/} (see
Definition~\ref{modulo} below). It is possible to obtain both upper and lower
bounds for this type of equivalence. The aim of this paper is to obtain the
upper bound. The lower bound has been found in~\cite{lowerbound}.

\begin{definition}\label{nest}\rm
A \emph{nest\/} of a vector field~$v$ (see Figure~\ref{fnest}) is an open
subset~$Z$ of its phase space such that
\begin{enumerate}
    \item $Z$ is homeomorphic to an annulus;
    \item the boundary curves of~$Z$ are limit cycles of the field;
    \item $Z$ contains no singular points of the field.
\end{enumerate}
\end{definition}

\begin{figure}
\begin{center}
\includegraphics{figure.1}
\end{center}
\caption{A nest.}\label{fnest}
\end{figure}

\begin{definition}\label{modulo}\rm
Let us consider vector fields~$v_1$ and~$v_2$. Let~$Z_1$ ($Z_2$) be the union
of all the nests of~$v_1$ ($v_2$). The fields~$v_1$ and~$v_2$ are called
\emph{equivalent modulo limit cycles\/} if the restriction of~$v_1$ to
$\R^2\setminus Z_1$ is orbitally topologically equivalent to the restriction
of~$v_2$ to $\R^2\setminus Z_2$.
\end{definition}

\subsection{Main Result}

\begin{theorem}\label{main1}
Denote by $K(H,n)$ the number of orbital topological equivalence types of
planar polynomial vector fields $v$ such that 1) $v$ has finite number of
singular points; 2) $v$ has at most $H$ limit cycles and 3) $\deg v\le n$.
Then
\[
K(H,n)\le C^{H+n^2},
\]
where $C=10^{157}$.
\end{theorem}

\begin{corollary}\label{hilbert}
For every $n>0$ the following statements are equivalent
\begin{enumerate}
\item The number of orbital topological equivalence classes of planar
polynomial vector fields with finite number of singular points and degree
less than or equal to~$n$ is finite. \item There is $H(n)$ such that every
planar polynomial vector field of degree less than or equal to~$n$ has at
most $H(n)$ limit cycles.

\end{enumerate}
\end{corollary}
\begin{proof}
$2\Longrightarrow1$. This is obvious from Theorem~\ref{main1}.

$1\Longrightarrow2$. We can restrict ourselves to $*$-fields. Indeed, if we
reduce components of a vector field by the common divisor, the number of
limit cycles can only increase (see the remark before
Definition~\ref{asterisk}).

Suppose that there are exactly~$m$ orbital topological classes of $*$-fields
of degree less than or equal to~$n$. Take one representative~$v_i$ for each
class ($i=1,\ldots,m$). Let~$H_i$ be the number of limit cycles of~$v_i$. By
the Finiteness Theorem (see~\cite{Finitiness1},~\cite{Finitiness2} and
also~\cite{History}) $H_i<\infty$. Thus
\[
H(n)=\max\limits_{1\le i\le m}H_i<\infty.
\]
\end{proof}

\begin{theorem}\label{main2}
Consider planar polynomial vector fields with finite number of singular
points and degree at most $n$. Denote by $M(n)$ the number of equivalence
classes modulo limit cycles of such fields. Then
\[
c^{n^2}\le M(n)\le C^{n^2},
\]
where $C=10^{471}$, $c=10^{10^{-8}}$.
\end{theorem}

Both Theorems above follow from a general Theorem below.

\begin{definition}\label{pfield}\rm
A closed set $\Pi\subset\R^2$ is called \emph{admissible\/} for a
$*$-field~$v$ if its boundary is a union of some limit cycles of~$v$. The
restriction~$w$ of the $*$-field to an admissible set~$\Pi$ is called a
\emph{$P$-field}. All the limit cycles of~$v$ that are components of
$\partial\Pi$ are considered limit cycles of~$w$. The \emph{degree\/} of $w$
is the degree of~$v$.
\end{definition}

\begin{theorem}\label{high}
Consider $P$-fields of degree at most~$n$ with at most~$H$ limit cycles.
There are at most
\[
C^{H+n^2}
\]
orbital topological equivalence classes of such fields, where $C=10^{157}$.
\end{theorem}

\begin{proof}[Proof of Theorem~\ref{main1}]
Notice that $\R^2$ is an admissible set for any $*$-field.
\end{proof}

\begin{proof}[Proof of Theorem~\ref{main2}]

The lower bound in the Theorem is proved in~\cite{lowerbound}. The explicit
constant $c=10^{10^{-8}}$ has not been written out but the calculation is
straightforward.

Take pairwise non-equivalent modulo limit cycles $*$-fields $v_1,\ldots,v_M$
with $\deg v_i\le n$ for $i=1,\ldots,M$. Let~$D_i$ be the union of all the
nests of~$v_i$. The set $\R^2\setminus D_i$ is admissible for~$v_i$, since
$D_i$ is the disjoint union of all the maximal nests of~$v_i$. Let~$w_i$ be
the restriction of~$v_i$ to $\R^2\setminus D_i$. By Definition~\ref{modulo}
the fields~$w_i$ are pairwise orbitally topologically non-equivalent.

\begin{lemma}\label{cyklicity}
Sum of the number of maximal nests of a $*$-field and the number of its limit
cycles that do not belong to any nest is less than or equal to the number of
singular points of this $*$-field.
\end{lemma}

It follows from this Lemma and from Bezout's Theorem that the number of limit
cycles of~$w_i$ is at most $2n^2$ (note that every maximal nest gives two
boundary limit cycles). Now we can apply Theorem~\ref{high}:
\[
M\le C^{2n^2+n^2}=\left(C^3\right)^{n^2},
\]
where $C=10^{157}$. Hence, $M(n)\le\left(C^3\right)^{n^2}=10^{471(H+n^2)}$.
\end{proof}

\begin{proof}[Proof of Lemma~\ref{cyklicity}]
Choose a limit cycle in every maximal nest. Call these limit cycles and the
limit cycles that do not belong to the nests \emph{labelled\/}. Our goal is
to assign a singular point to every labelled limit cycle.

Take any labelled cycle~$L$. Let $L_1,\ldots,L_s$ be all the labelled cycles,
satisfying two conditions:
\begin{enumerate}
\item $L_i$ is inside~$L$. \item There is no labelled cycle $L'$ such that
$L'$ is inside~$L$ and~$L_i$ is inside~$L'$.
\end{enumerate}
Let~$D$ be the domain bounded by~$L$, let~$D_i$ be the domain bounded by
$L_i$. Consider
\[
\Omega=D\setminus\overline{(D_1\cup\ldots\cup D_s)}.
\]

We claim that~$\Omega$ contains at least one singular point. Indeed, if
$s>1$, then it follows from the fact that Euler characteristic of $\Omega$ is
not equal to 0 (see Figure~\ref{lcycles}). Suppose $s=1$ and there are no
singular points in~$\Omega$. Then~$\Omega$ is a nest, so~$L$ and~$L_i$ are in
the same maximal nest. We come to contradiction.

\begin{figure}
\begin{center}
\includegraphics{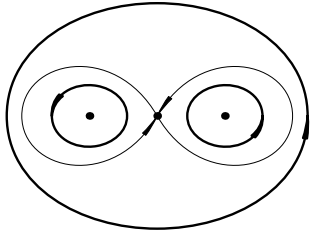}
\end{center}
\caption{Invariant set of nonzero Euler characteristic must contain a
singular point.}\label{lcycles}
\end{figure}

Assign to~$L$ any singular point inside~$\Omega$. Clearly, we have assigned
different singular points to different labelled cycles.
\end{proof}

The rest of the paper is devoted to the proof of Theorem~\ref{high}. The
paper is organized as follows: in Section~\ref{SingularPoints} we recall the
topological classification of singular points, define \emph{complexity\/} of
a singular point, and estimate the sum of complexities of all the singular
points of a $P$-field, using the method, which we have learnt
from~\cite{Khovanskii}. In Section~\ref{structure} we assign a graph on the
sphere to a $P$-field. It is proved in Appendix, that this graph is a
complete invariant of orbital topological classification of $P$-fields. The
proof is essentially using Theorems~75 and~76 of~\cite{ktds} (see
also~\cite{M1} and~\cite{M2}). In Section~\ref{combinatorics} we estimate the
number of graphs on the sphere, using the main result of~\cite{Tutte}.

The problem has been posed by A.~G.~Khovanskii and Yu.~S.~Ilyashenko, I want
to express my gratitude to them. Without the help of Yu.~S.~Ilyashenko this
paper would be unlikely completed. I would like to thank G.~Iyer for
correcting my English. I am grateful to the Department of Mechanics and
Mathematics of the Moscow State University, the Mathematical College of the
Independent University of Moscow, and the mathematical department of the
University of Chicago for their hospitality.


\section{Singular points}\label{SingularPoints}

A $P$-field cannot be flat at any point. Recall (see~\cite{ds1}, Chapter~5,
\S3) that every non-flat singular point~$O$ of a smooth planar vector field
is either monodromic or has a characteristic trajectory (i.e. a trajectory
that tends to~$O$ as $t\to+\infty$ or $t\to-\infty$, being tangent to some
line at~$O$).

\subsection{Classification of monodromic singular points}\label{monodromic}
Topological type of monodromic singular point is determined by its Poincare
map. This map cannot have infinite number of isolated fixed points due to the
Finiteness Theorem, since these fixed points correspond to limit cycles
(see~\cite{Finitiness1}, \cite{Finitiness2} and also~\cite{History}). Thus
every monodromic singular point of a $P$-field is either a focus, or a
center.

\subsection{Classification of characteristic singular points}\label{chrct}
\begin{figure}
\begin{center}
\includegraphics{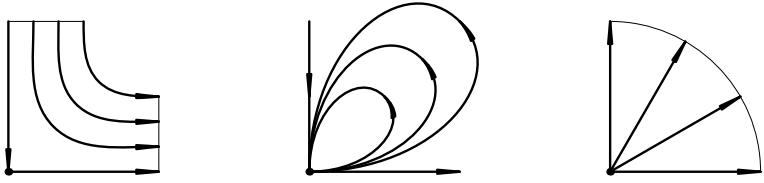}
\end{center}
\caption{Hyperbolic, elliptic, and parabolic sectors.}\label{sectors}
\end{figure}

A small neighbourhood of a characteristic singular point~$O$ splits into the
union of standard sectors. There are three types of standard sectors:
hyperbolic, elliptic, and parabolic  (see~\cite{ds1}, Chapter~5, \S3 and
Figure~\ref{sectors}).

\begin{remark}\rm
This splitting is not canonical: if we shrink the neighbourhood, then some
parts of the elliptic sector will become parabolic sectors. Nevertheless
Definition~\ref{complexity} and Definition~\ref{separatrix} are
``invariant''.
\end{remark}

\subsection{Compactification of the phase space}\label{sphere}

We want to compactify the phase space by adding an infinite point to the
plane. To get a chart in the neighbourhood of this point we identify the
plane with the 2-sphere, punctured at the North Pole, by stereographic
projection. The projection from the South Pole gives a chart in the
neighbourhood of the infinite point. The transition between these charts is
given by
\[
(x,y)\mapsto\left(\frac x{x^2+y^2},\frac{y}{x^2+y^2}\right).
\]
Thus in the second chart our vector field has the form
\[
\frac{v_1(x,y)}{(x^2+y^2)^n},
\]
where $v_1$ is a polynomial vector field. Set
\begin{equation}\label{infinity}
v_\infty(x,y)=(x^2+y^2)v_1(x,y).
\end{equation}
The fields of lines, corresponding to~$v$ and~$v_\infty$, agree on
$\R^2\setminus(0,0)$. Thus they can be patched together to the
field~$v_\mathrm{dir}$ on the sphere.

The infinite point is a singular point of~$v_\mathrm{dir}$, thus the
trajectories of~$v_\mathrm{dir}$ are those of~$v$ plus the infinite point.
This is why we need an extra $x^2+y^2$ factor in (\ref{infinity}) (otherwise
we could have a trajectory of~$v_\mathrm{dir}$ corresponding to two
trajectories of~$v$).

The following convention will be taken: the infinite point is considered a
singular point of~$v$.

\subsection{Complexity of singular points}

\begin{definition}\label{complexity}\rm
\emph{Complexity\/} of a characteristic singular point is the number of its
hyperbolic and elliptic sectors. Complexity of a monodromic singular point is
zero.
\end{definition}

We are going to estimate the total complexity (i.e. the sum of complexities)
of singular points. First we estimate the total complexity of all the
singular points except the point at infinity.

\begin{proposition}\label{compl}
The total complexity of finite singular points of a degree~$n$  $P$-field is
at most $6n^2-2n$.
\end{proposition}

\begin{proof}
Consider a $P$-field of degree $n\ge1$
\[
v(x,y)=P(x,y)\pa x+Q(x,y)\pa y.
\]

Set $R=P^2+Q^2$. Then $R(x,y)=0$ is the equation for the set of singular
points of~$v$.

\begin{lemma}\label{curve}
Let~$\Omega$ be any open bounded set of the plane, containing all the
singular points of~$v$. Let $\gamma_\epsilon$ be the curve given by the
equation $R(x,y)=\epsilon$. Then
\begin{enumerate}
\item For sufficiently small $\epsilon>0$ the curve $\gamma_\epsilon$  is
nonsingular in~$\Omega$.
\item The ovals of this curve tend to singular
points of~$v$ as $\epsilon\to0$.
\item For each singular point there is an
oval, enveloping this point.
\end{enumerate}
\end{lemma}

Consider the set~$S$ of points on~$\gamma_\epsilon$ where~$v$ is tangent to
this curve. This set is given by the system of equations
\[
\left\{
\begin{array}{l}
L_v(R-\epsilon)=0\\
R-\epsilon=0,
\end{array}
\right.
\]
where~$L_v$ is the directional derivative in the direction of~$v$. These
equations are polynomial of degree $3n-1$ and $2n$ respectively. Thus $S$ is
algebraic. Let~$\epsilon$ be so small that~$\gamma_\epsilon$ is nonsingular.
Clearly, if an oval of~$\gamma_\epsilon$ contains infinitely many points
of~$S$, then this oval is contained in~$S$. In this case the oval is a cycle
of~$v$. Thus~$S$ consists of isolated points and cycles of~$v$.

Now Bezout's Theorem tells us that~$S$ has at most $6n^2-2n$ isolated points.
It is an easy consequence of Lemma~\ref{curve} that in every hyperbolic or
elliptic sector of a singular point of~$v$ there is a point of
$\gamma_\epsilon$ where~$v$ is tangent to this curve; for $\epsilon$ small
enough these points are distinct and isolated (the latter follows from the
fact that a characteristic singular point has a neighbourhood without
cycles). Thus the overall number of elliptic and hyperbolic sectors is
bounded by $6n^2-2n$.
\end{proof}

\begin{proof}[Proof of Lemma~\ref{curve}]
The singular set of the curve~$\gamma_\epsilon$ is its intersection with the
set, given by the equation $\grad R=0$. We claim that every singular point
$O$ has a punctured neighbourhood without points where $\grad R=0$. If not,
then there is a smooth curve in this set such that~$O$ is in its closure,
since every algebraic set is a union of finite number of smooth strata.

Since~$R$ does not change along this curve, the curve consists of singular
points of~$v$ but this is impossible. Thus there is a punctured neighbourhood
of~$O$ without points of the set $\grad R=0$, and 1) is proved.

Let $O_1,\ldots,O_k$ be all the singular points of~$v$. Choose pairwise
disjoint open neighbourhoods~$U_i$ such that $\overline{U_i}\subset\Omega$.
Set
\[
\delta=\min\{R(x,y):(x,y)\in\overline\Omega\setminus\cup_{i=1}^k U_i\}.
\]
Clearly $\delta>0$. Take any $\epsilon<\delta$. We have
\[
\{(x,y):R(x,y)<\epsilon\}\cap\Omega\subset\cup_{i=1}^k U_i.
\]
This proves 2).

It remains to show that in~$U_i$ there is an oval,  enveloping~$O_i$. But if
there is no such an oval, then we can join~$O_i$ with some point of
$\Omega\setminus\cup_{i=1}^k U_i$ by a curve that does not intersect
$\gamma_\epsilon$, this contradicts continuity of~$R$.
\end{proof}

\begin{proposition}\label{compl2}
The complexity of the infinite point is at most $2n+2$.
\end{proposition}

\begin{proof}
It is similar to the proof of Proposition~\ref{compl}. The only difference is
that we should use the curve $x^2+y^2=C$ with~$C$ large enough instead of the
curve $P^2+Q^2=\epsilon$.
\end{proof}

\section{The structure of $P$-fields}\label{structure}

We shall assign an equipped oriented graph on the sphere  (possibly
disconnected) to every $P$-field. This graph can have loops and multiple
edges. This graph will be a complete invariant of the orbital topological
type of a $P$-field. Two graphs are considered equivalent if they are
isomorphic as equipped graphs, embedded into the sphere.

\begin{definition}\label{separatrix}\rm
A \emph{separatrix\/} of a $P$-field is a boundary trajectory of a hyperbolic
sector of a characteristic singular point.
\end{definition}

\subsection{Small graph}

We want to assign an edge of the graph to every separatrix of the $P$-field.
The problem is that the $\alpha$-limit set or the $\omega$-limit set of the
separatrix can consist of more than one point. Also, we want to have the
graph $C^1$-smooth. Thus we first introduce the notion of \emph{small graph}.

\begin{definition}\rm
Let $a\in\R$ or $a=-\infty$. An \emph{$\alpha$-germ\/} at~$a$ is an
equivalence class of maps $(a,b)\to S^2$, where two maps $(a,b_1)\to S^2$ and
$(a,b_2)\to S^2$ are equivalent if their restrictions onto $(a,b_3)$ coincide
for some~$b_3$. The \emph{$\omega$-germ\/} is defined similarly.

An $\alpha$-germ at $a\ne-\infty$ is called \emph{$C^1$-smooth\/} if some
(and then any) of its representatives can be extended to a $C^1$-smooth map
$(a-\epsilon,b)\to S^2$ for some $\epsilon>0$.

For an $\alpha$-germ at $-\infty$ we consider a representative
$\gamma:(-\infty,b)\to S^2$. The $\alpha$-germ is $C^1$-smooth if the map
$\gamma\circ\tan:(-\frac\pi2,\tan^{-1}b)\to S^2$ can be extended to a
$C^1$-smooth map $(-\frac\pi2-\epsilon,\tan^{-1}b)\to S^2$ for some
$\epsilon>0$. The similar definitions apply to $\omega$-germs.
\end{definition}

Let $\gamma$ be a trajectory of a $P$-field. Assume that $\gamma$ is neither
a singular point, nor a cycle, then it has a natural parametrization
$\gamma:(a,b)\to S^2$ (possibly $a=-\infty$, $b=+\infty$). The
\emph{$\alpha$-germ\/} of the trajectory is the $\alpha$-germ of~$\gamma$
at~$a$, the \emph{$\omega$-germ\/} is the $\omega$-germ of~$\gamma$ at~$b$.

\begin{definition}\rm
A separatrix of a $P$-field is called \emph{nice\/} if both its $\alpha$-germ
and $\omega$-germ are $C^1$-smooth. Other separatrices are called
\emph{nasty}.
\end{definition}

\begin{lemma}\label{nice}
The $\alpha$-germ of a trajectory~$\tau$ is $C^1$-smooth provided this
trajectory tends to a characteristic singular point~$O$ as $t\to-\infty$. The
similar statement is valid for $\omega$-germs.
\end{lemma}
\begin{proof}
By Theorem~64 of~\S20 of~\cite{ktds},~$\tau$ tends to~$O$ in a definite
direction. Choose an affine coordinate system in the neighbourhood of~$O$ so
that the direction is parallel to the $x$-axis. Suppose the $\alpha$-germ of
$\tau$ is not $C^1$-smooth.

\emph{Step 1.} We claim that there is a direction, not parallel to $x$-axis,
such that the points on~$\tau$ where the tangent line is parallel to this
direction accumulate to~$O$. Indeed, if in no neighbourhood of~$O$ the
projection of $\tau$ to $x$-axis is one-to-one, then we can take the
direction of $y$-axis. Otherwise in some neighbourhood of~$O$ the curve
$\tau$ is given by an equation $y=g(x)$, $x>0$. We have
\begin{equation}\label{deriv}
\frac{g(x)}x\to0,\mbox{ as }x\to0^+,
\end{equation}
since $\tau$ is tangent to $x$-axis. Since the $\alpha$-germ of~$\tau$ is not
$C^1$-smooth, $\lim\limits_{x\to0^+}g'(x)\ne0$. Indeed, otherwise we can
extend $g(x)$ to a $C^1$-smooth function, declaring $g(x)=0$ for $x\le0$. Now
it follows from (\ref{deriv}) and the Mean Value Theorem  that
$\lim\limits_{x\to0^+}g'(x)$ does not exist. Let $\lambda$ be any nonzero
number between the lower and upper limits of $g'(x)$. We can take direction
of the line $y=\lambda x$.

\emph{Step 2.} The set of points in the phase space where the vector field is
parallel to this direction is an algebraic set. Thus its intersection with
small enough neighbourhood of~$O$ is the union of smooth curves
$\gamma_1$,\dots,$\gamma_j$, where~$\gamma_i$ connects~$O$ with some
point~$O_i$. It is enough to show that for all~$i$ the intersection points of
$\gamma_i$ with~$\tau$ cannot accumulate to~$O$. This is clear if~$\gamma_i$
is not tangent to~$\tau$ at~$O$. Otherwise we shall use a version of
Rolle--Khovanskii method.

Assume that the points of intersection of $\tau$ with $\gamma_i$ accumulate
to~$O$. Since $\tau$ and $\gamma_i$ are analytic these points can accumulate
only to~$O$. Thus we can enumerate them in the order they appear on $\tau$:
\[
M_1, M_2, \ldots,
\]
where $M_s\to O$ as $s\to\infty$. If the field is tangent to~$\gamma_i$
everywhere, then $\gamma_i$ is part of a trajectory, and the claim is easy.
Otherwise, shrinking the neighbourhood of~$O$ we can assume that the field is
transversal to~$\gamma_i$ (except at $O$). Consider any curve~$\gamma_i'$,
satisfying the following properties: (1) It connects~$O$ with~$O_i$; (2) It
does not intersect~$\gamma_i$; (3) It is tangent to $y$-axis at~$O$.
Then~$\gamma_i$ and~$\gamma'_i$ bound a domain, denote it by~$\Omega$. The
part of~$\tau$ between~$M_s$ and~$M_{s+1}$ cannot be entirely in~$\Omega$,
since the field is transversal to~$\gamma_i$. Thus there is a point of~$\tau$
between~$M_s$ and~$M_{s+1}$, where~$\tau$ intersects~$\gamma'_i$.

We see that if the intersection points of~$\gamma_i$ with~$\tau$ accumulate
to~$O$, then so do the intersection points of~$\gamma'_i$ with~$\tau$. This
contradicts property (3) of~$\gamma'_i$.
\end{proof}

\begin{remark}\rm
The converse is also true but we do not need this.
\end{remark}

\begin{definition}\rm
Choose a point~$M_i$ on every limit cycle (including boundary limit cycles,
see Definition~\ref{pfield}). Choose one trajectory in every elliptic sector
of every singular point. The \emph{small graph\/} consists of the following
elements:
\begin{description}
    \item[Vertices:] The points~$M_i$ and the singular points of the $P$-field.
    \item[Edges:]
    The closures of the chosen trajectories in elliptic sectors, the limit
    cycles, and the closures of the nice separatrices.
\end{description}
\end{definition}

This graph is naturally embedded into the sphere. All the edges of the small
graph are $C^1$-smooth in this embedding (for trajectories in elliptic
sectors it follows from Lemma~\ref{nice}).

\subsection{Limit polycycles}

Consider a $C^1$-smooth curve on the 2-sphere. By \emph{side\/} of the curve
we mean the choice of co-orientation of this curve. Thus every $C^1$-smooth
curve has 2 sides.

\begin{definition}\rm
Let~$\gamma_1$ and~$\gamma_2$ be co-oriented separatrices of the vector field
(i.e. one of the two possible co-orientations on each separatrix is chosen).
We say that~$\gamma_2$ is the \emph{continuation\/} of~$\gamma_1$ if the
following holds:
\begin{enumerate}
    \item $\gamma_1$ tends to a singular point~$O$ as $t\to+\infty$,
    $\gamma_2$ tends to~$O$ as $t\to-\infty$;
    \item $\gamma_1$ and~$\gamma_2$ bound a hyperbolic sector of~$O$; this
    sector is on the positive side of each of the separatrices.
\end{enumerate}
\end{definition}

\begin{definition}\rm A \emph{limit polycycle\/} of a vector field is a cyclically
ordered finite set of singular points (with possible repetitions), and a
cyclically ordered set of disjoint co-oriented separatrices such that
\begin{enumerate}
    \item the time oriented $j$-th separatrix connects the $j$-th and $(j+1)$-st
    singular points;
    \item the $(j+1)$-st separatrix is the continuation of the $j$-th
    separatrix.
\end{enumerate}
(compare with~\cite{History}, \S3.1.)
\end{definition}

Consider a half-interval with the vertex on the limit polycycle (not at a
singular point). Assume that this half-interval is transversal to the field
everywhere and it is on the positive side of the polycycle. The monodromy
map~$g$ is defined in some neighbourhood of the vertex. Let~$z$ be the
coordinate on the half-interval ($z=0$ at the vertex, $z>0$ outside of the
vertex).

The fixed points of~$g$ correspond to limit cycles, thus there are only
finitely many of them by the Finiteness Theorem (see~\cite{Finitiness1},
\cite{Finitiness2} and also~\cite{History}). Thus for all~$z$ small enough we
have $g(z)>z$, $g(z)=z$ or $g(z)<z$. In these cases we call the limit
polycycle \emph{$\alpha$-limit, $0$-limit and $\omega$-limit\/} respectively.

Clearly, an $\alpha$-limit ($\omega$-limit) polycycle is the $\alpha$-limit
($\omega$-limit) set for all nearby trajectories that are on its positive
side. The converse is also true.

\begin{lemma}\label{limitset}
The $\alpha$-limit ($\omega$-limit) set of any trajectory is either a point,
or a limit cycle, or a limit polycycle.
\end{lemma}

\begin{proof}
Consider any infinite limit set~$X$ such that~$X$ is not a limit cycle. By
Theorem~68 of~\S23 of~\cite{ktds} the separatrices, constituting~$X$, can be
cyclically ordered so that for any two consecutive separatrices the next is
the continuation of the previous. The continuation is defined differently
in~\cite{ktds} (see Definitions~19 and~20 of~\S15). We leave to the reader to
check that this definition matches our definition. It follows that~$X$ is a
limit polycycle.
\end{proof}

\begin{figure}
\begin{center}
\includegraphics{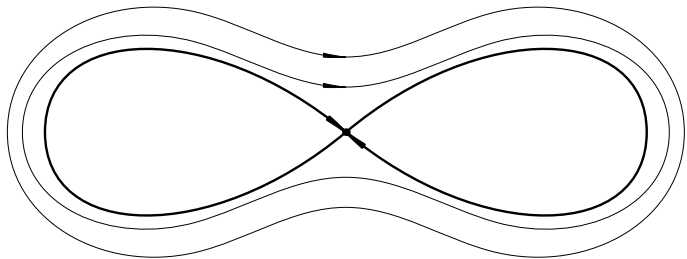}
\end{center}
\caption{A 0-limit polycycle.}\label{zerolimit}
\end{figure}

\subsection{Cycles without contact and large graphs}\label{without}

\begin{definition}\rm
A \emph{cycle without contact\/} of a vector field is a smooth closed curve
transversal to the vector field at every point.
\end{definition}

It is easy to see (see~\cite{ktds}, Lemma~2 of~\S24 and Lemma~1 of~\S 27)
that we can assign a cycle without contact to every $\alpha$-limit polycycle,
every $\omega$-limit polycycle, every focus, and every side of a limit cycle
so that
\begin{enumerate}
    \item Consider the region bounded by a limit polycycle, a focus, or a
    limit cycle and the corresponding cycle without contact. This region
    contains neither singular points, nor limit cycles, nor other chosen cycles
    without contact.
    \item Every trajectory, tending to a polycycle, to a focus, or to a limit
    cycle, intersects the corresponding cycle without contact at exactly one
    point. Every trajectory that intersects a cycle without contact tends to
    corresponding focus, polycycle, or a limit cycle.
    \item These cycles without contact intersect neither each other nor the
    elements of the small graph.
\end{enumerate}

Fix such a system of cycles without contact.

\begin{lemma}\label{nasty}
Every nasty separatrix intersects exactly one cycle without contact at
exactly one point.
\end{lemma}
\begin{proof}
We can assume that $\omega$-limit set of a nasty separatrix~$\gamma$ is a
characteristic limit point. It follows from Lemmas~\ref{nice}
and~\ref{limitset} that its $\alpha$-limit set is a focus, a limit polycycle,
or a limit cycle. Denote the corresponding cycle without contact by~$Y$.
Then~$\gamma$ intersects~$Y$ at exactly one point. It remains to prove that
$\gamma$ does not intersect other cycles without contact. Indeed, if $\gamma$
intersects another cycle~$Y'$, then the corresponding focus, limit cycle, or
limit polycycle should coincide with $\alpha$-limit set of $\gamma$.

Thus $Y$ and $Y'$ correspond to the same $\alpha$-limit set. There are only
two situations when it happens: (1) $Y$ and $Y'$ are the cycles without
contact corresponding to the same limit cycle and (2) $Y$ and $Y'$ correspond
to limit polycycles which coincide as sets but have different co-orientation
of trajectories. The first case is impossible, since $\gamma$ would have to
intersect this limit cycle. In the second case the polycycle is necessarily
homeomorphic to the circle; $\gamma$ would have to intersect the polycycle,
which is impossible.
\end{proof}
Thus a nasty separatrix is cut by the point of its intersection with the
corresponding cycle without contact into two parts. The part that is the
boundary trajectory of a hyperbolic sector is called \emph{truncated
separatrix}.

\begin{definition}\rm
The \emph{large graph\/} has all the vertices and edges of the small graph
and the following
\begin{description}
    \item[Vertices:]
    The intersection points of the nasty separatrices with the chosen cycles
without contact.
    \item[Edges:]
    The truncated nasty separatrices and the segments of the cycles without
contact in which they are split by the vertices, oriented counterclockwise
(if the cycle intersects no separatrix, then this cycle is not added to the
large graph).
\end{description}
\end{definition}

\begin{figure}
\begin{center}
\includegraphics{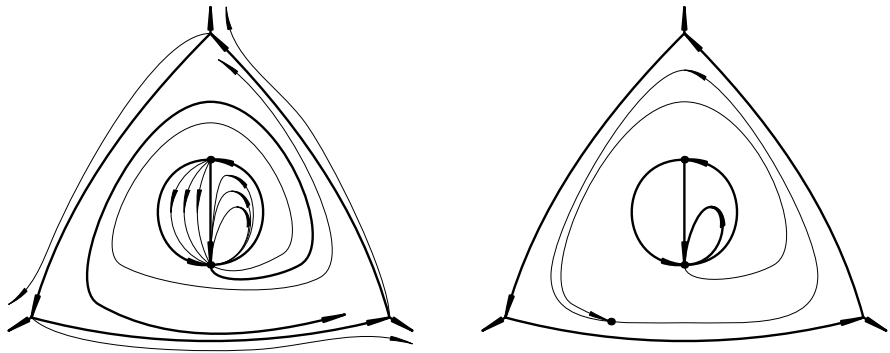}
\end{center}
\caption{A phase portrait and the corresponding large graph (thick curves on
the right picture correspond to the edges of the small graph).}
\end{figure}

\subsection{Equipping the graph}\label{equip}

Now we want to equip vertices and edges of the graph with some data. We
indicate on a vertex whether it is a repeller, an attractor, a clockwise
rotation center, a counterclockwise rotation center, or none of the above
(e.g. a saddle). Thus there are \textbf{5} possible equippings of a vertex.

On every edge we indicate whether it is a part of a cycle without contact, a
trajectory of an elliptic sector, a separatrix, or an edge, corresponding to
a limit cycle.

Consider a limit polycycle~$X$. Its separatrix~$\psi$ has a distinguished
side with respect to this polycycle. We mark this side of~$\psi$ by the
symbol~$\alpha$, $\omega$, or 0, depending on the type of polycycle.

It the side of a separatrix is not marked with~$\alpha$, $\omega$, or~$0$
(e.g. the separatrix does not belong to any limit polycycle) we mark this
side with $\emptyset$. In the same way we equip the edge, corresponding to a
limit cycle. However, only~$\alpha$, $\omega$, and~$\emptyset$ are possible
for a side of a limit cycle. Indeed, since the Poincare map of a limit cycle
is analytic there are only two possibilities: this map is either identity, or
the set of its fixed points is discrete. It cannot be identity, since the
cycle is limit, thus 0 is impossible on a limit cycle. The combination
$\emptyset\emptyset$ is also impossible. The other combinations are possible:
for example, $\omega\emptyset$ means that the limit cycle is a stable
boundary limit cycle (see Definition~\ref{pfield}). Thus there are \textbf{8}
possible equippings of a limit cycle.

For a separatrix all the combinations are possible. If a separatrix marked
with anything except $\emptyset\emptyset$, then it is nice and its
$\alpha$-germ and $\omega$-germ are both boundary trajectories of some
hyperbolic sectors by definition of limit polycycle.

If a separatrix is marked with $\emptyset\emptyset$, then we shall indicate
if its $\alpha$-germ and $\omega$-germ are boundary trajectories of
hyperbolic sectors. There are \textbf{3} cases because at least one germ
should be a boundary trajectory of a hyperbolic sector by definition of
separatrix. Notice that for a nasty separatrix the equipping is unnecessary
but we shall keep it for simplicity.

Thus there are $\mathbf{15+3=18}$ possible equippings of a separatrix. It
adds up to $\mathbf{1+1+8+18=28}$ possible equippings of an edge.

\begin{proposition}\label{technical}
Suppose that two $P$-fields have their equipped graphs isomorphic (as planar
oriented equipped graphs). Then these $P$-fields are orbitally topologically
equivalent.
\end{proposition}

The proof of this Proposition is a reformulation of the main result
of~\cite{ktds} and is given in the Appendix.

\section{Combinatorics}\label{combinatorics}

\begin{proposition}
Let $G_{conn}(l)$ be the number of connected graphs on the sphere (possibly
with loops and multiple edges) with at most~$l$ edges. Then
\begin{equation}\label{tutte}
G_{conn}(l)<12^l.
\end{equation}
(two graphs~$\Gamma_1$ and~$\Gamma_2$ are considered the same if they are
isomorphic as embedded graphs.)
\end{proposition}
\begin{proof}
The graph is said to be \emph{rooted\/} if one of its edges is chosen
together with its orientation and co-orientation (the other edges are neither
oriented, nor co-oriented). Let~$a_l$ be the number of rooted graphs with~$l$
edges. By~(5.1) of~\cite{Tutte}
\[
a_l=\frac{2(2l)!3^l}{l!(l+2)!}.
\]
A root is an additional structure. Thus if we forget the root, then we
decrease the number of graphs:
\[
G_{conn}(l)\le a_1+\ldots+a_l\le la_l=\frac{2l}{(l+2)(l+1)}\,3^l{2l\choose
l}<{12}^l.
\]
\end{proof}

\begin{remark}\rm
There is another approach due to physicists. Since every graph can be
completed to a triangulation, and the number of subgraphs of a triangulation
grows exponentially, it is enough to estimate the number of triangulations.
The following integral is the generating function for these numbers (for
graphs of different genera):
\[
\int_{\mathcal H_N}e^{-\frac tN\tr(H^3)}\,d\mu(H),
\]
where~$\mathcal H_N$ is the set of $N\times N$ hermitian matrices, $d\mu$ is
the Gaussian measure (see~\cite{Zvonkine}, \S7.1). However, the author has
never seen the proof of the required estimate, based on these methods.
\end{remark}

\begin{definition}\rm
The \emph{size\/} of a graph is the sum of the number of its edges and the
number of its vertices.
\end{definition}
\begin{proposition}\label{graphnum}
Let $G_{dc}(l)$ be the number of oriented possibly disconnected graphs on the
sphere (possibly with loops and multiple edges) with size at most~$l$. Then
\begin{equation}
G_{dc}(l)<48^l.
\end{equation}
\end{proposition}
\begin{proof}
Consider a graph~$\Gamma$ with~$b$ edges and~$f$ vertices ($b+f\le l$).
Suppose that this graph consists of~$d$ connected components, then $d\le f$.
Clearly, we can make the graph connected by adding $d-1$ edges. Denote the
new graph by $\Gamma'$. Then the number of edges of~$\Gamma'$ is less
than~$l$.

Thus every planar graph of size at most~$l$ can be obtained from a
\emph{connected\/} graph with at most~$l$ edges by deleting some edges. Since
there are at most $2^l$ subsets of edges of~$\Gamma'$, there are at most
$2^lG_{conn}(l)$ planar graphs with size~$l$ or less. Putting an orientation
on every edge, gives $2^b$ factor. Thus
\[
G_{dc}(l)\le2^b\cdot2^l\cdot G_{conn}(l)<{48}^l.
\]
\end{proof}

In Section~\ref{structure} we have assigned an equipped graph to every
$P$-field.
\begin{proposition}\label{size}
The size of the graph, corresponding to a $P$-field of degree~$n$ with~$H$
limit cycles, is at most $2H+37n^2+13$.
\end{proposition}
\begin{proof}
Denote by~$e$ and~$h$ the total numbers of elliptic and hyperbolic sectors
respectively. Denote by~$s_1$ and~$s_2$ the number of nice and nasty
separatrices respectively. Clearly, $s_1+s_2\le 2h$. The number of vertices
of the small graph is at most $H+(n^2+1)$ (remember the infinite singular
point!), the number of its edges is $H+s_1+e$.

The number of additional vertices of the large graph is~$s_2$, the number of
additional edges is~$2s_2$. Thus the size of the graph is at most
\begin{equation}\label{inequality}
6h+e+2H+n^2+1.
\end{equation}
By Propositions~\ref{compl} and~\ref{compl2} $e+h$ (this is the total
complexity of singular points) is at most $6n^2+2$ ($2n+2$ for the infinite
point, $6n^2-2n$ for the finite points). It remains to substitute the last
inequality in~(\ref{inequality}).
\end{proof}

Now we can finish the proof of Theorem~\ref{high}. The Theorem trivially
holds for $n=0$, since all the constant $P$-fields are topologically
equivalent. Thus we can assume $n>0$. Set $N=2H+37n^2+13$. Using
Propositions~\ref{graphnum} and~\ref{size}, we see that there are at most
$48^N$ possible graphs. Every vertex of the graph can have one of 5 possible
equippings, every edge of the graph can have one of 28 possible equippings,
thus there are at most $28^N$ possible equippings. Therefore there are at
most $(28\cdot48)^N$ possible equipped graphs.

Now we apply Proposition~\ref{technical}, it gives the estimate
$(28\cdot48)^N$ for the number of $P$-fields. Further,
\[
(28\cdot48)^N\le1344^{50(H+n^2)}<\left(10^{157}\right)^{H+n^2}.
\]
Theorem~\ref{high} is proved.


\appendix

\section{Appendix: Proof of Proposition~\ref{technical}}\label{ktds}

In~\cite{ktds} to every vector field on the sphere with finite number of
``singular elements'' a \emph{scheme\/} is assigned. The main Theorems
of~\cite{ktds} (Theorems~75 and~76 of \S29) assert that if two vector fields
have the same schemes, then they are orbitally topologically equivalent.
Unfortunately, the definition of scheme (Definition~33, \S29) is distributed
all over the book.

So our goal is to show that the scheme of a $P$-field can be recovered from
the large graph of the field. We do this examining subsequently all the
elements of the scheme. According to Definition~33 of \S29, we need to list
all the \emph{singular elements}, \emph{limit continua}, their \emph{global
schemes}, and all the pairs of \emph{conjugate free continua}. We shall
recall all the relevant definitions from~\cite{ktds}.

The other issue is that~\cite{ktds} deals with bounded phase spaces. However,
they mention that all the results are valid for a system on the sphere, see
\S29.5. We shall have to adjust some definitions of~\cite{ktds} to this case.

By default all the references in this Appendix are the references
to~\cite{ktds}.

\subsection{Singular elements}
According to Definition~33, there are~8 types of singular elements:

\textbf{1. Equilibrium states} are singular points in our terminology. They
are vertices of our graph. Notice that we can distinguish between singular
points and other vertices of the graph, using equippings of edges adjacent to
this vertices (if there are no such edges, then the vertex is necessarily a
singular point). Thus we can recover the list of singular points from the
graph.

\textbf{2. Orbitally unstable paths} A trajectory~$\gamma$ (it is called
\emph{path\/} in~\cite{ktds}) is called \emph{$\omega$-orbitally stable\/} at
a point $M\in\gamma$ provided that $\forall\epsilon>0\;\exists\delta>0$ such
that every trajectory~$\gamma'(t)$, passing through $\delta$-neighbourhood of
$M$ at $t=t_0$, remains in $\epsilon$-neighbourhood of~$\gamma$ for $t>t_0$,
\emph{$\alpha$-orbitally stable\/} trajectories are defined similarly. A
trajectory is called \emph{orbitally unstable} if it is not
$\alpha$-orbitally stable or not $\omega$-orbitally stable at least at one
point (see Definitions~14--17,~\S15). Note that in~\cite{ktds} this
definition is applied to bounded semitrajectories only. In order to make it
work on the sphere we have to use the spherical metric.

We claim that orbitally unstable trajectories are exactly limit cycles and
separatrices. Indeed, limit cycles are orbitally unstable by Theorem~37 of
\S15. It is quite clear that the separatrices are orbitally unstable.

Conversely, assume that~$\tau$ is an $\omega$-unstable trajectory that is not
a limit cycle. Theorem~40 of~\S15 shows that $\omega$-limit set of~$\tau$
consists of a single point. Theorem 38 of~\S15 tells that~$\tau$ is a
boundary curve of a hyperbolic sector. Thus $\tau$ is a separatrix. Therefore
we can recover the list of all orbitally unstable trajectories, using the
equippings of edges.

\textbf{3.} The remaining 6 types of singular elements deal with the boundary
of a region. Only so-called \emph{normal boundaries\/} are considered
in~\cite{ktds} (see~\S16.2). The normal boundary is one that consists of
finite number of arcs without contact and segments of trajectories (these
segments are called \emph{corner arcs}). These trajectories are not allowed
to be separatrices or limit cycles of the field.

Thus the boundary of a $P$-field is not normal. Therefore we have to do the
following: for a $P$-field~$v$ remove from the phase space the areas, bounded
by boundary limit cycles and their corresponding cycles without contact. Then
we get a new vector field~$v'$ with normal boundary. It is easy to see that
$P$-fields~$v_1$ and~$v_2$ are orbitally topologically equivalent if and only
if~$v_1'$ and~$v_2'$ are equivalent and the directions of rotation on the
corresponding boundary cycles of~$v_1$ and~$v_2$ are the same. Thus we shall
recover from graph the scheme of~$v'$ instead of that of~$v$.

\subsection{Scheme of a singular point}
The local scheme of a monodromic singular point can be trivially recovered
from the graph. The global scheme of a monodromic singular point is read from
the corresponding cycle without contact (see Proposition~\ref{corresp} of
this paper).

According to Definition~23 of~\S19, a scheme of a characteristic singular
point~$O$ is the list of \textbf{1)} all the separatrices of~$O$; \textbf{2)}
all the separatrices of other singular points that tend to~$O$; \textbf{3)}
elliptic sectors; \textbf{4)} the cyclic order of the above (recall that the
boundary does not have corner arcs).

Elliptic sectors correspond to the loops of the small graph. They can be
distinguished from other loops by their equipping. The separatrices are the
other edges, adjacent to~$O$. We can distinguish between separatrices of~$O$
and ``foreign'' separatrices, since we know for each germ of a separatrix
whether it is a boundary trajectory of a hyperbolic sector (see~\S\ref{equip}
of this paper). The cyclic order is specified, since the graph is embedded
into the sphere.

\subsection{Limit continua}

To comply with the terminology introduced in~\cite{ktds}, we use the
expressions \emph{$\alpha$-limit continuum\/} and \emph{$\omega$-limit
continuum\/} as the synonyms for $\alpha$-limit and $\omega$-limit sets of
trajectories.

A \emph{cell\/} of a vector field is a connected component of the phase space
after removal of all the singular elements. Consider the cell filled by
closed trajectories. This cell is doubly connected (see Theorem~50 of~\S16).
A connected component of its boundary is called \emph{$0$-limit continuum.}

One point limit continua are just attractors, repellers, and centers. Their
schemes can be read from the graph (easy). Hence, we shall restrict ourselves
to infinite limit continua.

\begin{proposition}\label{continua}
The infinite $\alpha$-limit, $\omega$-limit, and 0-limit continua that are
not limit cycles are limit polycycles and vice versa, limit polycycles are
limit continua. All the limit polycycles can be recovered from equipped
graph.
\end{proposition}
\begin{proof}
It follows from Lemma~\ref{limitset} of this paper that any infinite
$\omega$-limit continuum or $\alpha$-limit continuum is a limit cycle or a
limit polycycle. It can be proved similarly (using Theorem 70 of~\S23) that a
0-limit continuum is a 0-limit polycycle.

Conversely, it is clear that $\alpha$-limit and $\omega$-limit polycycles are
limit continua. Consider a 0-limit polycycle~$X$. We need to show that it is
a 0-limit continuum. Consider a half-interval where the monodromy map is
defined and the family of cycles, intersecting this half-interval. All these
cycles belong to the same cell; $X$ belongs to the boundary of this cell.
Thus~$X$ is a part of 0-limit continuum. It follows from Theorem~70 of~\S23
and the uniqueness of continuation of a separatrix that~$X$ coincide with
this 0-limit continuum.

It remains to show that the limit polycycles can be recovered from the graph.
A trajectory~$\psi$ belongs to some limit polycycle if and only if this
trajectory is equipped with~$0$, $\alpha$, or~$\omega$ on at least one of its
sides. It remains to show that we can ascertain from the graph whether two
separatrices belong to the same limit polycycle.

To this end we just need to check whether one separatrix is the continuation
of the other (because a co-oriented separatrix has at most one continuation).
This is the information we can get from the graph. Thus the whole limit
polycycle can be recovered from the graph.
\end{proof}

We know whether a continuum is $\alpha$-limit, $\omega$-limit, or $0$-limit
continuum from the equipping of any of its separatrices. The \emph{global
scheme \/} (see Definition 28, \S25) of the continuum is the list of all the
separatrices, tending to this continuum, with their cyclic order. This is
read from the graph by looking at the corresponding cycle without contact,
this is possible due to the following Proposition.
\begin{proposition}\label{corresp}
The correspondence between limit polycycles (limit cycles, foci) and cycles
without contact can be recovered from the graph.
\end{proposition}
\begin{proof}
The large graph splits the sphere into the parts, we shall call them
\emph{faces}. There is a natural embedding of the large graph into the phase
space. The part of the phase space, corresponding to the face under this
embedding, is called \emph{realization of a face}. Consider the case of
polycycle (the other cases are similar). The proposition follows from the
following claim: \emph{A limit polycycle~$X$ and a cycle without contact~$Y$
such that $Y$ belongs to the large graph, correspond to each other if and
only if
\begin{enumerate}
    \item They bound a face of the graph;
    \item This face is on the positive side of~$X$.
\end{enumerate}}

The ``only if'' statement follows from the choice of the cycles without
contact (see~\S\ref{without} of this paper).

Suppose~$X$ and~$Y$ bound a face and this face is on the positive side
of~$X$. Suppose, on the contrary, that~$Y$ corresponds to another limit
polycycle, limit cycle, or focus~$X'$. Since~$Y$ has been added to the graph,
there is a separatrix~$\tau$ such that $\tau$ intersects~$Y$. Then the
$\alpha$-limit set or the $\omega$-limit set of~$\tau$ coincides with~$X'$.
Consider the first case (the second case is similar).

Since~$\tau$ is a separatrix, its $\omega$-limit set must be a single point
$O$. But the $\omega$-limit set of~$\tau$ is contained in the closure of the
realization of a face, bounded by~$X$ and~$Y$, since~$Y$ is a cycle without
contact. Thus~$O$ belongs to~$X$ (hence, $X$ cannot be a limit cycle). But
this contradicts the assumption that~$X$ is a limit polycycle and the face is
on the positive side of~$X$.
\end{proof}

\subsection{Boundary Scheme}
The boundary (after removing neighbourhoods of boundary limit cycles)
consists of cycles without contact, labelled by boundary limit cycles (notice
that these cycles without contact may or may not belong to the large graph).

Let $Y$ be a cycle without contact, corresponding to a boundary limit
cycle~$X$. Since the boundary contains no corner arcs, the global scheme of
$Y$ is the list specifying (see Definition~30 of \S26): \textbf{1)}
whether~$Y$ is an outer or inner boundary curve; \textbf{2)} whether~$Y$ is
positive or negative cycle (i.e. whether $X$ is an $\alpha$-limit cycle or an
$\omega$-limit cycle); \textbf{3)} all singular paths, intersecting~$Y$,
enumerated in cyclic order.

Now \textbf{1)} and \textbf{2)} are read from the equipping of~$X$. The list
of all singular paths is recovered from vertices of the large graph (recall
that if~$Y$ does not belong to the large graph, then there are no singular
paths, intersecting~$Y$).

\subsection{Conjugate free continua}

Two 0-limit continua are called \emph{conjugate\/} if they bound a cell,
filled with closed trajectories. Two cycles without contact are called
conjugate if every trajectory that intersects the first one, intersects the
second one as well. An $\alpha$-limit continuum and an $\omega$-limit
continuum are called conjugate if the corresponding cycles without contact
are conjugate (see~\S28 and~\S27.4).

It remains to show that we can recover all the pairs of conjugate free
continua from the graph. (The continuum is called \emph{free\/} if no
separatrix intersects its cycle without contact. The conjugate continua are
obviously free.)

Again, consider the faces of the large graph. It is enough to prove the
following claim: \emph{limit continua~$X$ and~$X'$ are conjugate if and only
if they bound a face that is on the positive side of each continuum}. The
``only if'' part follows from Lemma~3a of \S28.2.

Conversely, suppose~$X$ and~$X'$ bound a face, and this face is on the
positive side of each of them. If the realization of this face is filled with
closed phase curves, then~$X$ and~$X'$ are conjugate 0-limit continua.
Otherwise let~$Y$ and~$Y'$ be cycles without contact, corresponding to~$X$
and~$X'$ respectively. We need to show that~$Y$ and~$Y'$ are conjugate
cycles. On the contrary, assume that~$\tau$ is a trajectory, $\tau$
intersects~$Y$ but does not intersect~$Y'$. Since~$Y$ is a cycle without
contact, $\tau$ intersects it just ones. Thus either positive, or negative
semitrajectory of~$\tau$ belongs entirely to the region bounded by~$Y$
and~$Y'$. But this region contains neither singular points nor limit cycles,
this contradicts the Poincare--Bendixon Theorem, which tells that
$\alpha$-limit ($\omega$-limit) set of any trajectory contains either a
singular point or a limit cycle.

\textsc{The University of Chicago, Department of Mathematics,\\
5734~S. University~Ave, Chicago, Illinois 60615, USA}

\emph{e-mail address}: \texttt{fedorov@mccme.ru}

\end{document}